\documentclass[a4paper,12pt]{article}

\newcommand{\proof}{\noindent{\bf Proof. }}
\newtheorem{thm}{Theorem}[section]   
\newtheorem{lem}[thm]{Lemma}
\newtheorem{cor}[thm]{Corollary}

\newcommand{\m}[4]
{{\footnotesize\pmatrix{
 #1 & #2 \cr
 #3 & #4 }}}

\title{Structure of Clifford Semigroups of Matrices\footnote{Research supported by
National Natural Foundation of China, grant 10571005. }}

\author{Yongwen Zhu}

\date{}

\begin{document}
\maketitle
\begin{center}

\noindent School of Mathematics and Information Science, Yanti University, Yantai City, P.R. China, 264005\\
E-mail: zyw@ytu.edu.cn 

\end{center}

\vskip12pt

 \footnotesize  {\bf Abstract:}\hskip6pt
 In this paper, we characterize completely the structure of Clifford semigroups of matrices over an arbitrary field. It is shown that a semigroups of matrices of finite order is a Clifford semigroup if and only if it is isomorphic to a subdirect product of some linear (0-)groups. Then we generalize this result to a semigroup of matrices of countably infinite order and give a similar necessary and sufficient condition for it to be a Clifford semigroup.

{\bf Key Words:}\hskip6pt
matrix semigroups; Clifford semigroups; subdirect product.
\par

{\bf MSC 2000:} \hskip6pt 20M25

\normalsize

\section{Introduction and preliminaries}

Throughout the paper, $F$ is an arbitrary field, $n$ is a positive integer. Let $M_n(F)$ and
$GL(n, F)$ be the set of all matrices and that of all invertible
matrices of $n\times n$ over $F$, respectively. A \emph{matrix semigroup}
is defined as a semigroup whose set is a subset of $M_n(F)$ and
whose composition is the usual multiplication of matrices, see \cite{okninski}, \cite{zhu-inverse}, \cite{zhu-compact}, \cite{zhu-regular}, \cite{zhu-complete simple} and \cite{zhu-Reducibility}.
If a matrix semigroup is a group, then we call it a \emph{matrix group}. For example,
$M_n(F)$ and $GL(n, F)$ are both matrix semigroups, and the latter
is even a matrix group. Evidently, any matrix semigroup must be a
subsemigroup of $M_n(F)$ . By a \emph{linear group}, we mean a subgroup of
$GL(n, F)$. Note that a matrix group must be
isomorphic to some linear group. For instance, $\{\m{x}{0}{0}{0}:x\in
F\setminus \{0\}
 \}$ is a group of matrices, which is isomorphic to $GL(1, F)$, the
linear group of dimension 1.

Inverse semigroups are an important class of semigroups, and Clifford semigroups are an important class of inverse semigroups, see \cite{howie}, \cite{petrich}, \cite{zhu-tiling} and \cite{maximal-clifford}. In fact, a \emph{Clifford semigroup} may be defined as a completely regular inverse semigroup.
We have studied the inverse semigroups of matrices in \cite{zhu-inverse}. In this paper, we shall characterize completely the structure of a Clifford semigroup of matrices. Our main result of this paper is Theorem \ref{main theorem}, which gives a sufficient and necessary condition for a semigroup of matrices of finite order to be a Clifford semigroup. This theorem should be generalized to those of matrices of infinite order in the last section of this paper.
As a direct corollary of our main theorem, Corollary \ref{maximal Clifford semigroup} improves the main results of \cite{maximal-clifford}.

For proving our main theorem in the next section, we present now some lemmas for preliminaries.

The first lemma comes from \cite[Corollary 3.5]{zhu-inverse}.

\begin{lem}\label{lemma about Clifford}
If $S$ is a Clifford semigroup of some $n\times n$ matrices with a
matrix of rank $r$, then there exists an invertible matrix $P$ so
that for all $a\in S$, there exist an $r\times r$ matrix $U$ and an
$(n-r)\times (n-r)$ matrix $Z$ such that $P^{-1}aP=\m{U}{0}{0}{Z}$.
\end{lem}

We call a matrix $a\in M_n(F)$ \emph{fullranked} if the rank of $a$ is $n$, that is $r(a)=n$. We call a semigroup $S\subseteq M_n(F)$ \emph{fullranked} if every element of $S$ is fullranked. If $S$ contains the zero matrix 0 and every nonzero matrix of $S$ is fullranked, then we say $S$ is \emph{0-fullranked}. If $G$ is a group, then $G^0=G\cup \{0\}$ is a semigroup. We call a segmiroup formed in this way a \emph{$0$-group}. Now we can prove our second lemma.

\begin{lem}\label{fullranked}
A (0-)fullranked regular semigroup of matrices is a (0-)group.
\end{lem}

\proof Suppose that $S\subseteq M_n(F)$ is fullranked regular semigroup, that is $S\subseteq GL_n(F)$. For any nonzero element $a\in S$, $a$ is actually an invertible matrix. So $aba=a$ implies that $ab=1$, the identity matrix. Therefore, $S$ is a linear group.

Similarly, one can prove that if $S\subseteq M_n(F)$ is 0-fullranked regular semigroup, then $S$ is a 0-group.
So the proof is complete ¨…

If $S\subseteq M_n(F)$ is a semigroup and $r(a)=r(b)$ for any $a,b\in S$, then we call $S$ \emph{monoranked}.
If $S\subseteq M_n(F)$ is a semigroup containing the zero matrix 0 and $r(a)=r(b)$ for any $a,b\in S\setminus \{0\}$, then we call $S$ \emph{0-monoranked}.
For Green's relations ${\cal H}$,${\cal R}$,${\cal L}$,${\cal D}$ and ${\cal J}$, we refer to \cite{howie}, \cite{petrich} and \cite{zhu-tiling}. Then we can show the next statement, the last lemma of this section.

\begin{lem}\label{monoranked}
A (0-)monoranked Clifford semigroup of matrices is a (0-)group.
\end{lem}

\proof Suppose that $S\subseteq M_n(F)$ is 0-monoranked Clifford semigroup. Then for any $a,b\in S\setminus \{0\}$, $r(a)=r(b)$. In light of \cite[Lemma 2.1]{okninski}, $a{\cal J}b$. From \cite[Theorem 4.2.1]{howie}, it follows that ${\cal J}={\cal H}$ for Clifford semigroup$S$. Thus $a{\cal H}b$. Therefore, $S\setminus \{0\}$ is a group and then $S$ is a 0-group.\par
Similarly, one can deduce that $S$ is a group provided that $S\subseteq M_n(F)$ is monoranked Clifford semigroup.
So the proof is complete ¨…

\section{Structure of Clifford semigroups of matrices of finite order}

Suppose that $S\subseteq M_n(F)$ is a semigroup and there exists a matrix $P\in GL(n,F)$ such that for every $a\in S$, $P^{-1}aP=diag\{a_1,\cdots, a_t,0,\cdots,0\}$ with $a_k \in M_{n_k}(F)$ for $k=1,\cdots,t$. For any $k\leq t$, let $S_k=\{a_k|a\in S\}$. Then $S_k\subseteq M_{n_k}(F)$ is a semigroup. We denote $S\sim [S_1,\cdots,S_k]$. Now we can prove our main theorem of this paper, which characterizes completely the Clifford semigroups of
$n\times n$ matrices, where $n$ is finite.

\begin{thm}\label{main theorem}
Let $F$ be a field and $S\subseteq M_n(F)$ a semigroup. Then $S$ is a Clifford semigroup if and only if $S$ is isomorphic to a subdirect product of some linear (0-)groups.
\end{thm}

\proof Sufficiency. If $S$ is isomorphic to a direct product of some (0-)groups, then $S$ is obviously a Clifford semigroup because every (0-)group is evidently a Clifford semigroup.\par
Necessity. Assume that $S\subseteq M_n(F)$ is a Clifford semigroup. If there is a matrix $a\in S$ with $0<r(a)<n$, then by Lemma \ref{lemma about Clifford}, we have $S\sim [S_1,S_2]$. The same argument can be used to $S_1$ and $S_2$ respectively. So by using repeatedly Lemma \ref{lemma about Clifford}, we should assume that $S\sim [S_1,\cdots, S_t]$ with each $S_i$ (0-)fullranked. It is easy to see that $S$ is Clifford semigroup if and only if $S_1,\cdots, S_t$ are all Clifford semigroups. Thus $S_1,\cdots, S_t$ are all (0-)fullranked regular semigroups of matrices. In light of Lemma \ref{fullranked}, each $S_i$ is a group or a 0-group. Therefore, $S$ is isomorphic to the subdirect product of linear (0-)groups $S_1,\cdots, S_t$ and the proof is complete. ¨…

The next corollary is a direct consequence of the last theorem, which appears also as \cite[Lemma 4.1]{zhu-inverse}.

\begin{cor}
A Clifford semigroup of $2\times 2$ matrices over a field $F$ is
commutative if it contains a matrix of rank 1.
\end{cor}

\begin{cor}
If $S\subseteq M_3(F)$  is a Clifford semigroup which contains at least one singular matrix, then every matrix of $S$ has at least one eigenvalue in $F$.
\end{cor}

\begin{cor}\label{rank 1 finite}
If $S\subseteq M_n(F)$ is a Clifford semigroup and every matrix of $S$ is of rank 1, then $S$ is a commutative group and all eigenvalues of every matrix of $S$ are in $F$.
\end{cor}

\proof Since $S$ is monoraked, $S$ is a group by Lemma \ref{monoranked}. According to Theorem \ref{main theorem}, one can easily complete the proof. ¨…

The next corollary is an immediate consequence of Theorem \ref{main theorem}.

\begin{cor}\label{maximal Clifford semigroup}
Let $F$ be a field and $S\subseteq M_n(F)$ a semigroup. Then $S$ is a maximal Clifford subsemigroup of $M_n(F)$ if and only if $S$ is isomorphic to the $0-$group $GL_n(F)\bigcup \{0\}$.
\end{cor}

Note that the last corollary generalizes Theorem 7 of \cite{maximal-clifford}, which stated the same result only under the condition that $F$ is a finite field.

\section{Structure of Clifford semigroups of matrices of countably infinite order}

Given a finite positive integer $n$, we define a function $\phi _n: M_n(F)\rightarrow M_{n+1}(F)$ as follows: $\phi _n (a)=\m{a}{0}{0}{1}$ for any $a\in M_n(F)$. It is clear that $\phi _n$ is a homomorphism of semigroups. Let $M_{{\cal X}_0}(F)=\overrightarrow{\lim }(M_n(F),\phi _n)$, the direct limit of the system $(M_n(F),\phi _n)$. Then  $M_{{\cal X}_0}(F)$ is a semigroup, called the \emph{full matrix semigroup of countably infinite order}.

For any $s\in M_{{\cal X}_0}(F)$, we can view $s$ naturally as a matrix with the size ${\cal X}_0\times {\cal X}_0$.

If $S\subseteq M_{{\cal X}_0}(F)$ is a semigroup, then we call $S$ is \emph{matrix semigroup of countably infinite order}. Now, we can generalize our main result presented in last section to semigroups of matrices of countably infinite order and obtain the following theorem.

\begin{thm}\label{generalization of main theorem}
Let $F$ be a field and $S\subseteq M_{{\cal X}_0}(F)$ a semigroup. Then $S$ is a Clifford semigroup if and only if $S$ is isomorphic to a subdirect product of some linear (0-)groups of finite order.
\end{thm}

This theorem can be proved by the same method as used in the proof of Theorem \ref{main theorem}, but the Lemma \ref{lemma about Clifford} should be replaced by the next lemma.

\begin{lem}\label{generalzation of lemma about Clifford}
If  $S\subseteq M_{{\cal X}_0}(F)$ is a Clifford semigroup with a
matrix of rank $r$, then there exists an invertible matrix $P$ of countably infinite order so
that for all $a\in S$, there exist an $r\times r$ matrix $U$ and a matrix $Z$ of countably infinite
such that $P^{-1}aP=\m{U}{0}{0}{Z}$.

\end{lem}

\proof Since $S$ is a Clifford semigroup, $S$ is an inverse semigroup. We denote its idempotents
by $E$. If $S$ contains a matrix $b$ of rank $r$, then $bb^{-1}=e\in
E$. It is clear that $r(e)=r(bb^{-1})=r(b)=r$.
By the definition of $M_{{\cal X}_0}(F)$, there exists a positive integer $n$ such that $e=\m{f}{0}{0}{0}$ with $f\in M_n(F)$. It follows that $f$ is an idempotent from that $e$ is an idempotent.
Evidently, $r(f)=r(e)=r$. Thus there exists $P_0\in GL(n, F)$ so that $P_0^{-1}fP_0=\m{I_r}{0}{0}{0}$. Let $P=\m{P_0}{0}{0}{I}$ where $I$ is a matrix of countably infinite order all elements of which are 0 but the diagonals are 1s. Then $P$ is an invertible matrix of countably infinite order, and $P^{-1}eP=\m{I_r}{0}{0}{0}$.

For any $a\in S$, set $P^{-1}aP=c$. Since $P^{-1}SP$ is isomorphic
to $S$ and $S$ is a Clifford semigroup, $P^{-1}SP$ is also a
Clifford semigroup. In light of \cite[Theorem 4.2.1]{howie}, in a Clifford
semigroup, any element and any idempotent are commutative, so
$ce=ec$. Assume $c=\m{U}{X}{Y}{Z}$ with $U\in M_r(F)$. Then
$\m{U}{X}{Y}{Z}\m{I_r}{0}{0}{0}=\m{I_r}{0}{0}{0}\m{U}{X}{Y}{Z}$,
that is $\m{U}{0}{Y}{0}=\m{U}{X}{0}{0}$. It follows that $X=0$ and $
Y=0$. So $P^{-1}aP=c=\m{U}{0}{0}{Z}$ as desired. ¨…

As an immediate consequence of Theorem \ref{generalization of main theorem}, the following corollary is similar to the Corollary \ref{rank 1 finite}.

\begin{cor}\label{rank 1 infinite}
If $S\subseteq M_{{\cal X}_0}(F)$ is a Clifford semigroup and every matrix of $S$ is of rank 1, then $S$ is a commutative group and all matrices of $S$ are diagonalizable simultaneously.
\end{cor}

{\bf Acknowledgement}£ºThe author would like to express his
gratitude to Professor Wang Efang.


\end{document}